\documentclass[11pt]{article}

\usepackage{amssymb}
\usepackage{latexsym}
\def\barr{\begin{array}}
\def\earr{\end{array}}
\def\ds{\displaystyle}

\def\reals{{{\rm l}\kern - .15em {\rm R} }}
\def\scriptO{{{\it O}\kern -.42em {\it `}\kern + .20em}}

\def\reals{{{\rm l}\kern - .15em {\rm R} }}
\def\caL{\mathcal{L}}
\def\gg{\mathfrak{g}}
\def\gh{\mathfrak{h}}
\def\ga{\mathfrak{a}}
\def\cc{\mathbb{C}}
\def\pp{\mathbb{P}}
\def\rr{\mathbb{R}}
\def\bp{\mbox{\bf P}}
\def\bc{\mbox{\bf C}}
\newtheorem{Thm}{Theorem}[section]

\newtheorem{Cor}[Thm]{Corollary}
\newtheorem{SubThm}{Theorem}[subsection]
\newtheorem{SubProp}[SubThm]{Proposition}
\newtheorem{Prop}{Proposition}[section]

\newtheorem{Lem}[Thm]{Lemma}
\begin{document}

\thispagestyle{empty}

\null
\medskip

\centerline{\Large \bf PROJECTIVE RANKS OF COMPACT}
\bigskip
\centerline{\Large \bf HERMITIAN SYMMETRIC SPACES}

\large
\vspace{1in}

\centerline{\bf Amassa Fauntleroy \hspace{-.08truein}}
\centerline{Department of Mathematics \hspace{-.08truein}}
\centerline{North Carolina State University \hspace{-.08truein}}
\centerline{Raleigh, NC 27695 \hspace{-.08truein}}
\centerline{U.S.A. \hspace{-.08truein}}

\vspace{1in}
\normalsize

\baselineskip20pt



\centerline{\bf Abstract}

Let $M$ be a compact irreducible Hermitian symmetric space and write $M=G/K$, with $G$ the group of holomorphic isometries of $M$ and $K$ the stability group of the point of $0 \in M$.  We determine the maximal dimension of a complex projective space embedded in $M$ as a totally geodesic submanifold.

\vspace{1.25in}


\pagebreak

\noindent{\Large \bf Introduction}
\medskip

Let $M$ be a simply connected compact complex manifold carrying an Hermitian metric of everywhere nonnegative holomorphic bisectional curvature.  In [11] Mok proved that if the second Betti number $b_2(M)=1$, then $M$ is biholomorphic to an irreducible compact Hermitan symmetric space.  It was earlier proved by Siu and Yau [14] and Mori [12] that when the above curvature is everywhere positive $M$ is biholomorphic to a complex projective space.  In each of the three papers cited above the existence of a minimally embedded projective line in $M$ plays a crucial role.

It turns out that such minimally embedded projective lines are totally geodesic in $M$.  In this paper we study the maximal totally geodesic complex submanifolds of $M$ which are biholomorphic to a complex projective space.  We call the dimension of such a submanifold the projective rank of $M$.  We calculate the projective ranks of each of the irreducible compact Hermitian symmetric spaces.  The results are given in Section 5.  In Sections 1 through 4 we develop the techniques used to make these calculations and in the last section we discuss the degrees of the holomorphic totally geodesic maps $\varphi : \mathbb{P}^r_{\mathbb{C}} \rightarrow M$ where $r = $ projective rank of $M$.  We also discuss the conjugacy of these maximal totally geodesic complex projective spaces in $M$ under the action of the group of isometries of $M$.

We shall make use of the work of Chen and Nagano [5] on totally geodesic submanifolds of symmetric spaces.  Several authors have studied the question of minimal or energy minimizing maps from $S^2$ to a compact symmetric space.  The interested reader may consult [2] and [3] for connections with the present paper.  It is the recent work of Robert Bryant [18] which motivates the author to revive these results which were summarized without proofs in the article [19].  Unlike Bryant, the approach taken here is very much an algebraic one with only an occasional nod to the topological and analytic methods which underlie many of the foundational results.  We thank the referee of an earlier version of this paper for pointing out an error in our original discussion of the case of the quadrics.

\section{Preliminaries on Hermitian Symmetric Spaces.}

A complex manifold $M_0$ (noncompact) with Hermitian metric $h$ is an Hermitian symmetric space or H.S.S. if each point of $M_0$ is an isolated fixed point of an involutive holomorphic isometry of $M_0$.  Let $G_0$ denote the connected component of the group of holomorphic isometries of $M_0$.  Then $G_0$ is a connected Lie group which acts transitively on $M_0$.  Fix $p_0 \in M_0$ and let $K \subset G_0$ be the isotropy group of $p_0$ so $M_0 = G_0/K$.  We introduce the following notations:

\bigskip

$\gg_0$: Lie algebra of $G_0$

$s_0$ : involutive isometry having $p_0$ as isolated fixed point

$\tilde{k}$ : Lie algebra of $K$

$\gg_0 = k + m_0$ : decomposition of $\gg_0$ into the $+1$ and $-1$ eigenspaces if $\sigma = ad(s_0)$

$\gg = \gg^{\mathbb{C}}_0$ : complexification of $\gg_0$

$m = m^{\mathbb{C}}_0$ : complexification of $m_0$

$\gg_c = k + m_c$ : compact real form of $\gg$ where $m_c = im_0$(cf., [7; V 2.1]).

\bigskip

Let $t$ be a Cartan subalgebra of $k$.  Then $t$ is a Cartan subalgebra of $\gg_0$ and $\gg_c$ and $\mathfrak{h} = t^{\mathbb{C}}$ is a Cartan subalgebra of $\gg$.  In fact, $G_0$ is a connected centerless semi-simple Lie group.  The center $Z(K)$ of $K$ is a torus whose dimension is the number of simple direct factors of $G_0$ and $K$ is the centralizer in $G_0$ of $Z(K)$. 

Let $\Phi$ be the $\gh$-root system of $\gg$ so

\bigskip

${\ds \gg =\gh \oplus \coprod_{\alpha \in \Phi} \gg_{\alpha}}$

$\Phi_K$ : compact roots; i.e., $\gh$-root system of $k^{\mathbb{C}}$

$\Phi_M$ : noncompact roots, i.e., $\{\alpha \in \Phi : \gg_\alpha \subseteq m \}$

$z$ : central element of $k$ such that $J$=ad($z$) induces the complex structure of $M_0$

$m^\pm$ : $(\pm i)$-eigenspace of $J$ on $m$

$p : k^{\mathbb{C}} + m^-$, parabolic subspace of $\gg$ that is the sum of nonnegative eigenspaces of ad($iz$) : $\gg \rightarrow \gg$

$P$: parabolic subgroup of $G$; i.e., the complex analytic subgroup corresponding to $p$.

The following theorem summarizes the major classical results which we will need.  Details and proofs may be found in [7], [17].

\begin{Thm}\quad Let the notation be as above.  Then 
\begin{description}
\item{(1)} If $b_2(M_0) = 1$ then $M_0$ is irreducible, $G_0$ is a compact centerless simple Lie group and $K$ is a maximal compact proper subgroup of $G_0$ and $G_c$ and has a 1-dimensional center.
\item{(2)} $P$ is a maximal parabolic subgroup of $G, G_c$ acts transitively on $M = G/P$ and $G_c \cap P = K$.
\item{(3)} $M_0 = G_0/K$ embeds holomorphically in $M$ as an open $G_0$-orbit, and $M$ is a compact H.S.S.
\item{(4)} Assume that $b_2(M) = 1$.  Let $\Delta \subset \Phi$ be a base for the root system $\Phi$.  Then there is an $\alpha \in \Delta$ such that
\[
	P=P(\Delta - \{\alpha\}) \qquad (cf.,[10; Section \ 1])
\]
Moreover, the homogeneous line bundle $L_\alpha$ on $G/P$ corresponding to $\alpha$ is ample.  Thus $G/P$ is projective algebraic.
\end{description}
\end{Thm}

Let $Y\subset M$ be a submanifold containing $p_0$.  Then $Y$ is totally geodesic in $M$ if and only if $T_{p_0 Y} = n$ satisfies $[[n,n],n] \subset n$ -- i.e., $n$ is a Lie triple system.  In this case $n+[n,n]$ is a Lie subalgebra of $\gg_c$.  If $H$ is the corresponding analytic subgroup of $G_c$ then $H$ acts transitively on $Y$ and $Y$ is again a symmetric space.  Note:  $Y$ need not be Hermitian symmetric.

\begin{Lem} \quad Let $X$ be a simply connected symmetric space and $Y \subset X$ a totally geodesic simply connected submanifold of $X$.  Assume that the connected group of isometries of $Y, \ I_0(Y)$, is semisimple.  Then every isometry of $Y$ extends to an isometry of $X$.
\end{Lem}

\noindent{\bf Proof.}\quad Fix $y \in Y$ and let $n=T_{y, Y}$ be the corresponding Lie triple system.  Let $\ga = n + [n,n]$ be the corresponding subalgebra of $\gg_X$ -- the Lie algebra of $I_0(X)$ where $I_0(X)$ is the connected component of the isometry group of $X$.  Let $A \subset I_0(X)$ be the analytic subgroup corresponding to $\ga$ and $K_A = A \cap K, \ K$ the stability group in $I_0(X)$ of $y$.  Then by [7, p. 225] $Y = A/K_A$.  In particular, there is a natural mapping $\beta : A \rightarrow I(Y)$ and by [7, Remark 2, p. 211] the image is a closed analytic subgroup of $I(Y)$.

Let $K'$ be the stability group of $y$ in $I(Y)$.  If $\gamma \in I(Y)$ and $\gamma \cdot y = y'$, then since $A$ acts transitively on $Y$, there is an $h \in A$ which $hy'=y$.  Thus $\beta(h) \gamma \in K'$.  Since $I_0(Y)$ is semisimple, the Lie algebra $k'$ of $K'$ is spanned by the set
\[
	\{R_0(u,v): u,v \in T_{y,Y} \}
\]
by [9; Vol. II, XI 3.2(5)].  But for $u,v$ in $T_{y,Y}$ we can write $u = d \beta (\overline{u}), v=d \beta(\overline{v})$ with $\overline{u}, \overline{v}$ in $\ga$.  Since on $X, R_y(u,v)(w) = -[[u,v], w]$ and $Y$ is totally geodesic, we have
\[
\barr{ll}
	R_0(u,v)(w) &= -[[u,v], w] \\
		      &=-[[d \beta(\overline{u}), d \beta (\overline{v})], \; w]
\earr
\]
If follows that $d \beta [\ga, \ga] = k'$.  Hence
\[
	d \beta : \ga \rightarrow \ \mbox{Lie}(I_0(Y))
\]
is surjective.  Since the image of $\beta$ is closed, $\beta(A) = I_0(Y)$.

Now $I_0(Y)$ is simply connected to there exists a subgroup $A_0 \subseteq A$ such that $\beta: A_0 \rightarrow I_0(Y)$ is an isomorphism and $I_0(Y) \subseteq I(X)$.  This proves the lemma.

\begin{Cor} \quad Let $Y$ be a totally geodesic Hermitian symmetric subspace of the H.S.S. M.  Assume that $Y$ and $M$ are compact and that $Y$ is irreducible.  Let $G' = A(Y)$ and $G=A(M)$ denote the respective linear groups of holomorphic transformations of $Y$ and $M$.  Then there exists an isogeny $\tilde{G}' \rightarrow G'$ and a homomorphism $\tilde{G}' \rightarrow G$ such that 
\[	
	Y = \tilde{G}' / \tilde{P}' \rightarrow M=G/P
\]
is the inclusion map.
\end{Cor}

\noindent{\bf Proof.} \quad Since the connected groups of holomorphic isometries $I_0(Y)$ and $I_0(M)$ are dense in the connected components of $G'$ and $G$ (cf., [7, p. 211]) respectively, the corollary follows readily from the lemma.

\noindent{\bf Remark.} \quad Note that if $Y$ is irreducible and $\gg = k + m$ is the Cartan decomposition of the Lie algebra of $I_0(Y)$ then $m$ generates $\gg$, i.e., $\gg = m + [m,m]$.

We give now an alternative description of Hermitian symmetric spaces in terms of $C$-spaces.  Let $G$ be a connected complex simply connected Lie group with simple Lie algebra $\gg$ of rank $\ell$.  Let $\gh$ be a Cartan subalgebra of $\gg$ and let $\Phi$ denote the set of roots of $\gg$ so
\[
	\gg = \gh + \sum_{\alpha \in \Phi} \gg_\alpha
\]
Fix a base $\Delta = \{ \alpha_1, \ldots, \alpha_\ell \}$ for $\Phi$ and let $\Phi^+$ denote the set of positive roots. Fix a simple root $\alpha_r, \ 1 \leq r \leq \ell$ and put
\[ \barr{rl}
	\Phi_1 &= \{ \alpha \in \Phi : n_r(\alpha) = 0\} \\
\noalign{\medskip}
	\Phi(n^+) & = \{ \alpha \in \Phi^+ : n_r (\alpha) > 0 \} \\
\noalign{\medskip}
	\Phi(n) & = \Phi_1 \cup \Phi (n^+) 
\earr \]
Define subalgebras of $\gg$ by
\[
\barr{rl}
	\gg_1 & = \gh + {\ds \sum_{\alpha \in \Phi_1}} \gg_\alpha \\
\noalign{\medskip}
	n^+ & = {\ds \sum_{\alpha \in \Phi (n^+)}} \gg_\alpha \\
\noalign{\medskip}
	n & = \gh + {\ds \sum_{\alpha \in \Phi(n)}} \gg_2
\earr
\]
Then $\gg_1$ is reductive, $n^+$ is nilpotent and $n = \gg_1 + n^+$.  Let $\gg_c$ be a compact real form of $\gg$ with $\gg_{\mathbb{R}} \cap \sqrt{-1} \gh_{\mathbb{R}} \neq (0)$.  Let $P$ be the complex analytic subgroup of $G$ such that Lie $P = n$ and $G_c$ the real analytic subgroup with Lie $G_c = \gg_c$.  Put $K = G_c \cap P$.  Then $G_c$ acts transitively on $G/P$ and induces the structure of a compact complex manifold on the homogenous space $G_c/K \simeq G/P$.  The irreducible compact Hermitian symmetric spaces are given as follows:

(A III) \quad Grassmannians : $Gr(r, \ell +1) = (A_\ell, \alpha_r)$

(BDI) \quad	Quadric Hypersurfaces:
\[
	Q^N = \left\{ \barr{l}
	(B_\ell, \alpha_1) \ \mbox{if} \ N= 2 \ell -1 \\
	(D_\ell, \alpha_1) \ \mbox{if} \ N = 2 \ell -2 
	\earr \right.
\]

(CI) \quad $(C_\ell, \alpha_\ell)$

D III \quad $(D_\ell, \alpha_\ell)$

E III \quad $(E_6, \alpha_1)$

E VII \quad $(E_7, \alpha_7)$

Note that the group $P$ defined above is a maximal parabolic subgroup -- the standard parabolic defined by the subset $\Delta_r = \Delta - \{ \alpha_r \}$ of $\Delta$.  Moreover, $K$ is a maximal compact subgroup of $G_c$ and the Levi-factor $L$ of $P$ is the centralizer in $G$ of $T_r = \mbox{Ker} \ \alpha_r$.  The Lie algebra of $L$ is $k + ik, \ k = \mbox{Lie} \ K$.

A somewhat more concrete description of the classical compact H.S.S. is as follows: (cf., [17, p. 321].

(A III) \quad Grassmannians:  $G_c = SU (\ell + 1), K = SU( \ell + 1) \cap \{ U(r) \times U(\ell + 1 - r) \}, \ G_c/K$ is the space of $r$-dimensional subspaces of $ \mbox{\bf C}^{\ell + 1}$

(BDI) \quad The nonsingular quadratic hypersurface in $\bp^{N + 1}$ defined by
\[
	z^2_0 + \cdots + z^2_{N+1} = 0
\]

$G_c = SO(N +2), K = SO(N) \times SO(2)$

(CI) \quad $Sp(n)/U(n) \subset G(n, 2n)$

The space consists of $n$-dimensional linear subspaces of $ \mbox{\bf C}^{2n}$ annihilated by a nondegenerate skew-symmetric bilinear form.

(D III) \quad The subvariety $G_c/K = SO(2n)/U(n)$ of the same Grassmannian $G(n, 2n)$ consisting of $n-$dimensional subspaces annihilated by a nondegenerate symmetric bilinear form.

In the sequel we sometimes want to distinguish between $Gr(d,n)$ thought of as $d$-dimensional linear subspaces of $\bp^n_{\mathbb{C}}$ and as $d+1$ dimensional subspaces of $\mbox{\bf C}^{n+1}$.  In those instances we use the notation $G(s, V)$ or $G(s,n)$ for the Grassmannian thought of as $s$-dimensional subspaces of the $n$-dimensional vector space V.

\section{Minimal and Linear Embeddings of $\bp^r_{\mathbb{C}}$ in Complex Grassmannians.}

Let $V$ be an $n+1$ dimensional complex vector space and $d$ a positive integer with $n \leq 2d < 2n$.  Let $Gr(d,n)$ denote the Grassmann manifold of $d+1$ dimensional subspaces of $V$ or equivalently $d$-dimensional linear subspaces of $\bp(V) \simeq \bp^n_{\mathbb{C}}$.  The manifold $Gr(d,n)$ is Hermitian symmetric of type A III.  If $Y$ is a compact complex submanifold of $Gr(d,n)$ which has everywhere positive holomorphic bisectional curvature then by [12] or [14] $Y$ is biholomorphic to a complex projective space.  If $Y$ is also totally geodesic in $Gr(d,n)$ then $Y$ is a symmetric subspace of $Gr(d,n)$.  In this section we study such submanifolds $Y$ which have minimal degree in $Gr(d,n)$.

Let $E(d,n)$ denote the universal subbundle over $Gr(d,n)$ and $Q(d,n)$ the universal quotient bundle.  Then
\[
0 \rightarrow E(d,n) \rightarrow Gr(d,n) \times V \rightarrow Q(d,n) \rightarrow 0
\]
is exact and the fiber over $p \in Gr(d,n)$ of $E(d,n) \rightarrow Gr(d,n)$ is the subspace of $V$ represented by the point $p$.  For a holomorphic embedding $\varphi : \bp^r_{\mathbb{C}} \rightarrow Gr(d,n)$ we define the \underline{degree of $\varphi$} to be the degree of the line bundle $\det(\varphi^* E^v(d,n))$ where $E^v(d,n)$ is the dual of $E(d,n)$.  As usual for a sequence of integers $0 \leq a_0 < a_1 < \ldots < a_d \leq n$ we define the Schubert variety $\Omega(a_0, a_1, \ldots, a_d)$ as follows: Fix a sequence of subspaces
\[
	V_0 \subset V_1 \subset \cdots \subset V_d \subset P(V)
\]
with $\dim V_r = a_r$.  Then
\[
	\Omega(a_0, \cdots, a_d) = \{L \in Gr(d,n)| \dim(L \cap V_i) \geq i \}
\]
The classes of these cycles in $H_*(Gr(d,n), \mbox{\bf Z})$ generate the homology.  Suppose $a_i = i +1$ for $0 \leq i \leq d$.    Then the Schubert variety $\Omega(1,2,\ldots, d,d+1)$ is just the set of subspaces $L \subset V$ which are contained in $V_d$.  Hence $\Omega(1,2,\ldots, d + 1) = Gr(d, V_d) = \bp^{d+1}$. 

The Pl\"{u}cker embedding of $Gr(d,n)$ is determined by the line bundle $\det(Q(d,n)) = \det(E^V(d,n))$.  For any Schubert variety $\Omega(a_0, \ldots, a_d)$ of dimension $k$ we have (cf., [6; p. 274])
\[
	\deg \Omega(a_0, \cdots, a_d) = \frac{k!}{a_0!a_1! \cdots a_d!} \cdot \prod_{i<j} (a_j-a_i)
\]
where $k = {\ds \sum^{d}_{i=0}} a_i - \left(\barr{l} d \\ 2 \earr \right)$.  In particular for $\Omega(1,2, \ldots, d+1)$, $k = d+1$ and $\deg \Omega(1,2, \ldots, d+1) = 1$.

Now let $\varphi : \bp^r \rightarrow Gr(d,n)$ be a holomorphic embedding.  We say that $\varphi$ is \underline{minimal} if $\varphi^* c_1(Q(d,n))=[O_{\mathbb{P}r}(1)]$ and that $\varphi$ is \underline{linear} if $\varphi$ factors as $\varphi_0 \circ \varphi_1$ where $\varphi_1 : \bp^r \rightarrow \bp^{d+1}$ is a linear embedding of $\bp^r$ as a subspace of $\bp^{d+1}$ and $\varphi_0$ is an isomorphism $\bp^{d+1}$ with a Schubert variety $\Omega(1,2, \cdots, d+1)$ in $Gr(d,n)$.  The embedding $\varphi$ is \underline{totally geodesic} if $\varphi(\bp^r)$ is a totally geodesic submanifold of $Gr(d,n)$.

\begin{Lem} \quad Let $\varphi : \bp^1 \rightarrow Gr(d,n)$ be an embedding.  Then $\varphi$ is minimal if and only if it is linear.  
\end{Lem}

\noindent{\bf Proof.}  \quad A linear embedding is clearly minimal.  Conversely, if $\varphi$ is minimal then $\varphi^*Q$ is a rank $n - d$ vector bundle generated by its global sections.  By [13; p. 22],
\[
	\varphi^*Q = \mathcal{O}_{\mathbb{P}^1} (e_1) \oplus \cdots \oplus \mathcal{O}_{\mathbb{P}^1}
	(e_{n-d})
\]
with $e_1 \geq e_2 \geq \cdots e_{n-d} \geq 0$.  Since $\det \varphi^* Q = \mathcal{O}_{\mathbb{P}^1}(1)$ we have $e_1=1$ and $e_i =0$ for $i>1$.  Now consider the pull back of the universal sequence
\[
	0 \rightarrow E(d,n) \rightarrow (Gr(d,n) \times \bc^{n+1}
	\stackrel{\pi}{\rightarrow}  Q(d,n) \rightarrow 0
\]
to $\bp^1$.  We have
\[
	0 \rightarrow E(d,n) | _{\mathbb{P}^1} \rightarrow \bc^{n+1} \times \bp^1 \rightarrow \mathcal{O}_{\pp^1}(1) \oplus 
	\mathcal{O}_{\mathbb{P}}^{n-d-1} \rightarrow 0
\]
Let $\{s_0, \ldots, s_n\}$ be a global frame for the trivial bundle $\bc^{n+1} \times Gr(d,n)$ such that $\pi(s_{d+2}, \ldots, \pi(s_n)$ span the trivial $\mathcal{O}^{n-d-1}$ of $\varphi^*Q$ over $\bp^1$.  For any point $x \in \bp^1$ the  corresponding subspace $L_x \subset \bc^{n+1} \times \bp^1$ must lie in the subspace $B$ spanned by $\langle s_0, \ldots, s_{d+1} \rangle$.  Thus $\varphi( \bp^1) \subset Gr(d,B)$ hence $P$ is linear.  \hfill $\Box$

For a line $\mathcal{L}$ in $Gr(d,d+1) \simeq \bp^{d+1}$ we have a simple description of $\mathcal{L}$ as a Schubert variety.  Let $B=C^{d+2}$.  There is a $d$-dimensional subspace $B_{d-1}$ of $B$ such that
\[
\barr{ll}
	\mathcal{L} & = \{L \subset B | B_{d-1} \subset L \} \\
\noalign{\medskip}
	& = \Omega (0,1,2, \ldots, d,d+2)
\earr
\]
In particular, if $L_1, L_2$ are in $\mathcal{L}$ then
\[
	\dim(L_1 +L_2) \leq d +1 +d +1 -d =d +2
\]
So if $L_1 \neq L_2, \ L_1 + L_2 = B$

\begin{Lem} \quad Let $\varphi : \bp^2 \rightarrow Gr(d,n)$ be a minimal embedding.  Then $\varphi$ is linear.
\end{Lem}

\noindent{\bf Proof.}  \quad Let $\mathcal{L}, \mathcal{L}'$ denote two distinct lines in $\bp^2$.  Then there exists subspaces $B$ and $B'$ of $\bc^{n+1}$ of dimension $d+2$ and subspaces $B_{d-1} \subset B, \ B'_{d-1} \subset B'$ of dimension $d$ such that
\[
\barr{ll}
	\varphi(\mathcal{L}) & = \{L|B_{d-1} \subset L \subset B \} \\
\noalign{\medskip}
	\varphi(\mathcal{L}') & = \{ L' | B'_{d-1} \subset L' \subset B' \}
\earr
\]
I claim $B = B'$.  Suppose not.  Then since $\caL \cap \caL' \neq \emptyset$ there is a common subspace $L_0$ of $B$ and $B'$ of dimension $d+1$.  Thus $\dim(B + B') = d+2+d+2-(d+1)=d+3$ and $\varphi(\caL)$, $\varphi(\caL')$ both lie in $G(d+1, B'') = Gr(d,P(B''))$ where $B''=B+B'$.  Since $B_{d-1}$ and $B'_{d-1}$ are both subspaces of $L_0$ we can find a basis $\{e_0, \ldots, e_{d+2} \}$ of $B''$ such that the following conditions are fulfilled
\[	
\barr{rl}
	B_{d-1} &= \langle e_0, \ldots, e_{d-1} \rangle \\
\noalign{\medskip}
	B'_{d-1} &= \langle e_1, \ldots, e_d \rangle \\
\noalign{\medskip}
	L_0 & = \langle e_0, \ldots, e_d \rangle \\
\noalign{\medskip}
	B &= L_0 + \langle e_{d+1} \rangle \\
\noalign{\medskip}
	B' &= L_0 + \langle e_{d+2} \rangle
\earr
\]
Then the subspace $L_1 = \langle e_0, \ldots, e_{d-1}, e_{d+1} \rangle$ is in $\caL$ and $L_2 = \langle e_1, \ldots, e_{d+2} \rangle$ is in $\caL'$.  Hence
\[
	\dim(L_1 + L_2) = d +3
\]
But $L_1$ and $L_2$ correspond to points in $\varphi(\bp^2)$ so lie in a minimally embedded line $\caL''$.  By Lemma 2.1 $\caL''$ is linearly embedded so by the above remark
\[	
	\dim(L_1 + L_2) = d + 2 
\]
This contradiction leads to the desired conclusion $B=B'$ and hence $\varphi(\bp^2) \subset Gr(d + 1, B')$ so is linearly embedded. \hfill$\Box$

Let us fix a basis $e_0, \ldots, e_n$ of $V$ establishing a fixed isomorphism with $\bc^{n+1}$.  Let $B \subset V$ be the subspace spanned by $e_0, \ldots, e_{d+1}$.  Then without loss of generality one may assume that a linearly embedded $\bp^{d+1}$ has image $G(d+1, B)$.  Restricting the universal sequence to $\bp^{d+1}$ gives
\[
	0 \rightarrow E(d,n) |_{\mathbb{P}^{d+1}} \rightarrow \bc^{n+1} \times \bp^{d+1} \rightarrow
	Q(d,n)|_{\mathbb{P}^{d+1}} \rightarrow 0
\]
If $x \in \bp^{d+1}$ then the corresponding subspace $L_x$ lies in $B$ so the $n-d-1$ vectors $e_{d+2}, \ldots, e_n$ thought of as global sections of $\bc^{n+1} \times \bp^{d+1}$ map onto global sections of $Q(d,n)|_{\mathbb{P}^{d+1}}$ free from relations.  This implies that
\[
	Q|_{\mathbb{P}^{d+1}} = L \oplus (\bc^{n-d-1} \times \bp^{d+1})
\]
where $L$ is a line bundle.  Thus since $c_1(\varphi^*Q(d,n))=1, \ \varphi^*Q(d,n) = \mathcal{O}_{\mathbb{P}^{d+1}} (1) \oplus \mathcal{O}^{(n-d-1)}_{\mathbb{P}^{d+1}}$.  Conversely, if $\varphi^* Q(d,n) = \mathcal{O}_{\mathbb{P}^{d+1}} (1) \oplus \mathcal{O}^{(n-d-1)}_{\mathbb{P}^{d+1}}$ we can find $e_{d+2}, \ldots, e_n$ global sections of $\bc^{n+1} \times Gr(d,n)$ which map onto generators of the trivial factor of $Q(d,n)|_{\mathbb{P}^{d+1}}$.  Extending these to a basis $e_0, \ldots, e_d$ for $\bc^{n+1} \times Gr(d,n)$ yields $L_x \subset \langle e_0, \ldots, e_{d+1} \rangle$ for each $x \in \bp^{d+1}$.

\begin{Thm} \quad Let $\varphi : \bp^r \rightarrow Gr(d,n)$ be an embedding with $1 \leq r \leq d + 1$.  Then $\varphi$ is minimal if and only if it is linear.
\end{Thm}

\noindent{\bf Proof.} \quad It suffices by the above remark to show that $\varphi^*Q(d,n) = \mathcal{O}_{\mathbb{P}^{d+1}} (1) \oplus \mathcal{O}^{(n-d-1)}_{\mathbb{P}^{d+1}}$.  Now for $r=1$ or $2$ the Lemmas 2.1, 2.2 yield the result.  If $2 < r < d + 1$ we have for any line $\caL \subset \bp^r$, $\varphi/\caL$ is minimal so $\varphi^*Q|_{\caL} = \mathcal{O}_{\caL} (1) \oplus \mathcal{O}^{(n-d-1)}_{\caL}$.  It follows that $\varphi^* Q(d,n)$ is uniform - i.e., for each $\caL \subset \bp^r \ \varphi^* Q/\caL = \mathcal{O}_{\caL} (a_1) \oplus \cdots \oplus \mathcal{O}_{\caL}(a_{n-d})$ with $a_i$ constant.  By Lemma 2.2, $\varphi^*Q(d,n)$ splits on a $\bp^2 \subset \bp^r$ also.  Then by a theorem of Horrocks [13, p. 42] $\varphi^* Q(d,n)$ splits over $\bp^r$ so $\varphi$ is linear. \hfill $\Box$

The above arguments actually yield a more general result.

\begin{Cor} \quad Let $Q$ be a vector bundle on $\bp^n_{\mathbb{C}}$ of rank $d$ and generated by $n+d$ global sections.  If $Q$ is uniform of type $(1,0, \ldots, 0)$ then $Q$ splits on $\bp^n$, i.e., $Q \simeq \mathcal{O}_{\mathbb{P}}^n (1) \oplus \mathcal{O}^{(d-1)}_{\mathbb{P}^n}$.
\end{Cor}

\noindent{\bf Proof.} By hypothesis there exist a surjection
\[
	\mathcal{O}^{n+d}_{\mathbb{P}^n} \stackrel{\varphi}{\longrightarrow} Q \rightarrow 0
\]
Let (*): $0 \rightarrow E \rightarrow \mathcal{O}^{n+d}_{\mathbb{P}^n} \stackrel{\varphi}{\longrightarrow} Q \rightarrow 0$ be the corresponding exact sequence of locally free sheaves.  Over $\bp^1 \subset \bp^n$ we obtain
\[
	0 \rightarrow E|_{\mathbb{P}^1} \rightarrow \mathcal{O}^{n+d}_{\mathbb{P}^1} \rightarrow
	\mathcal{O}_{\mathbb{P}^1} (1) \oplus \mathcal{O}^{d-1}_{\mathbb{P}^1} \rightarrow 0 \;.
\]
It follows that the morphism $\psi : \bp^n \rightarrow Gr(n,n +d)$ defined by (*) restricts to a linear embedding on each $\bp^1 \subset \bp^n$. The proof of Lemma 2.2 shows that $\psi$ is then linear on each $\bp^2 \subset \bp^n$.  By the remarks preceding Theorem 2.3 $Q = \psi^*Q(n,d+n)$ splits on $\bp^2$ into $\mathcal{O}_{\mathbb{P}}(1) \oplus \mathcal{O}^{d-1}_{\mathbb{P}^2}$.  By the theorem of Horrocks [13, p. 42] $Q$ splits over $\bp^n$. \hfill $\Box$

In addition, we can obtain the following restriction on the existence of maps from projective spaces into Grassmannians.

\begin{Cor} \quad If $K>d+1$ then every holomorphic map $\varphi$ from $\bp^K_{\mathbb{C}}$ to $Gr(d,n)$ with $c_1(\varphi^*Q(d,n))=1$ is constant.
\end{Cor}

\noindent{\bf Proof.} Suppose $K>d+1$ and let $W \subset \bp^K_{\mathbb{C}}$ be a linear subspace of dimension $d+2$.  If $\varphi$ is holomorphic map from $\bp^K$ to $Gr(d,n)$ with $c_1(\varphi^*Q(d,n)) = 1$ then for every hyperplane $W' \subseteq W$, the restriction of $\varphi$ to $W'$ is a minimal hence linear embedding.  Suppose $W_1$ and $W_2$ are two distinct hyperplanes in $W$ and let $L_i$ be a line in $W_i$ such that $L_1 \cap L_2 \neq \emptyset$ and $L_1 \not\subset W_2, L_2 \not\subset W_1$.  We can find $B_i \subset \bc^{n+1}, \  i =1,2$ such that
\[
	\varphi(W_i) = G(d+1, B_i)
\]
and $B'_i \subset B_i$ of dimension $d$ such that
\[
	\varphi (L_i) = \{L \in G(d+1, B_i) | B'_i \subset L \subset B_i\} \;.
\]
Just as in the proof of Lemma 2.2 we can conclude that $B_1 = B_2 = B$.  It follows that $\varphi(W') \subset G(d+1, B) \subset Gr(d,n)$ for every hyperplane $W' \subset W$ so $\varphi (W) \subset Gr(d+1, B)$.  But then $\dim W=d+2$ and the only morphism from $W \cong \bp^{d+2}$ to $G(d+1, B) \simeq \bp^{d+1}$ are known to be constant.  \hfill $\Box$

The interested reader may compare the above result with the results of Tango in [15] and [16].

\section{Some Representation Theory.}

Let $X$ be a compact H.S.S. and $Y$ a totally geodesic compact Hermitian symmetric subspace of $X$.  If $Y$ is isomorphic to $\bp^\ell_{\mathbb{C}}$ as H.S.S., then we have seen in Corollary 1.3 that there is a homomorphism
\[
	\rho : SU(\ell +1, \bc) \rightarrow G
\]
where $G$ is the group of analytic automorphisms of $X$.  If $X$ is a Grassmannian, say $Gr(d,n)$, then $\rho$ yields a homomorphism which we again call $\rho$
\[
	\rho : G_Y = SL(\ell + 1, \bc) \rightarrow SL(n+1, \bc)
\]
so $\rho$ is a representation of $G_Y$ on $V = \bc^{n+1}$. We want to describe these representations when $2 \ell \geq n$.  

Let $\lambda_1, \cdots, \lambda_\ell$ be the fundamental dominant weights for $SL(\ell + 1, \bc)$.  We fix the maximal torus $T \subset SL(\ell + 1, \bc)$ to be the diagonal subgroup and use the standard description of the roots as in [8, p.64].  In particular, if $d= \mbox{diag}(d_1, \ldots, d_n) \in T$ the base for the root system (B = upper triangular matrices) is given by
\[
	\alpha_i(d) = d_i d^{-1}_{i+1} \quad 1 \leq i \leq n
\]
and
\[
	\lambda_i = \frac{1}{\ell + 1} [(\ell - i + 1) \alpha_1 + \cdots + (i-1)(\ell -i+1) \alpha_{i-1} 
	+ i(\ell - i+1) \alpha_i + i(\ell - i) \alpha_{i+1} + \cdots + i \alpha_\ell ]
\]
for $1 \leq i \leq \ell$.  A straightforward computation gives
\[	
	\lambda_i(d) = d_1 d_2 \cdots d_i \:.
\]
If $e_1, \ldots, e_{\ell+1}$ are the standard basis elements for $V=C^{n+1}$ then $e_1 \wedge \cdots \wedge e_i \in \stackrel{i}{\Lambda} V$ is the highest weight vector corresponding to $\lambda_i$.  Hence
\[
	\deg(\lambda_i) = \left( \barr{c} n+1 \\ i \earr \right)
\]
Recall that for a dominant weight $\lambda = m_1 \lambda_1 + \cdots + m_\ell \lambda_\ell$ Weyl's character formula (cf. [7, VII, 10.2, 11.2]) gives
\[
	\deg \lambda = \left\{ \prod_{\alpha \in \Phi^+} \langle \lambda + \delta, \alpha \rangle \right\} \bigg/
	\prod_{\alpha \in \phi^+} \langle \delta, \alpha \rangle \;.
\]
Following Humphries [ibid] we can use the co-root $\alpha^{\vee}$ instead of $\alpha$.  Then for $\alpha$ height $r$
\[
	(\lambda + \delta, \alpha^{\vee}) = m_i + \cdots + m_{i+r-1} + r
\]
if $\alpha = \alpha_i + a_{i+1} + \cdots + \alpha_{i +(r-1)}$ is a root of height $r$ and

\noindent $(*) \hspace{1.5in} \deg \lambda = \left\{ {\ds \prod_{\alpha \in \Phi^+}} (\lambda + \delta, \alpha^{\vee}) \right\} \bigg/ {\ds \prod_{\alpha \in \Phi^+}} (\delta, \alpha^{\vee})$

\begin{Prop} \quad Let $G=SL(\ell + 1, \bc)$ with $\ell \geq 4$ and let $V$ be a $G$-module of dimension $d$.  Assume that $\ell + 1 < d < 2 \ell$.  Then $V$ contains a trivial $G$-module of dimension $r \geq d - \ell - 1$ as a direct summand.  If $\ell \geq 5$ and $\ell + 1 < d < 2(\ell + 1)$ then $V$ contains a trivial summand of dimension $r \geq d - \ell - 1$.
\end{Prop}

\noindent{\bf Proof.} Let $W \subset V$ be a nontrivial irreducible submodule of $V$ with highest weight $\lambda = m_1 \lambda_1 + \cdots + m_\ell \lambda_\ell$.  I claim that $m_1 + \cdots + m_\ell = 1$.  Since $(\lambda +  \delta, \alpha^V) = m_i + \cdots + m_{i+r-1} + r$ if $\alpha = \alpha_i + \cdots + \alpha_{i+r-1}$ it follows that if $\lambda > \lambda'$ in the lexicographic ordering then $\deg \lambda > \deg \lambda'$.  If $\lambda = {\ds \sum^{\ell}_{i=1}} m_i \lambda_i$ with $m_i >0$ for some index $i$ with $2 \leq i \leq \ell - 1$ then $\deg \lambda_i = \left(\barr{c} \ell + 1 \\ i \earr \right)$.  But $\left( \barr{c} \ell + 1 \\ i \earr \right) > 2(\ell +1 )$ for $i$ in this range and $\ell \geq 4$.  It follows that $m_i  = 0$ if $2 \leq i \leq \ell - 1$.

Consider now $\lambda = m_1 \lambda_1 + m_\ell \lambda_\ell$.  Again $\deg D > \min(\deg m_1 \lambda_1, \deg m_\ell \lambda_\ell)$.  If $m_1 \geq 2$ then
\[
	\deg(m_1 \lambda_1) = m_1 \cdot 1^{\ell - 1} \cdot (m_1 + 2)2^{\ell -3} \cdots
	(m_1 + \ell) / \Delta > m_1 \deg \lambda_1 \geq 2 \deg \lambda_1 = 2(\ell + 1)
\]
Similarly, if $m_\ell \geq 2$ then $\deg(m_\ell \lambda_\ell \geq 2 \deg \lambda_\ell = 2( \ell +1)$.  Thus the hypothesis $d \leq 2(\ell + 1)$ implies $m_1 \leq 1$ and $m_2 \leq 1$.  If $\lambda = \lambda_1 + \lambda_2$ then
\[
\barr{ll}
	\deg D = & \frac{2 \cdot 2^{\ell-2} \cdot 3^2 2^{\ell -3} \cdots (\ell +2)}{\Delta} \\
\noalign{\medskip}
	& >2 \deg \lambda_1 = 2 \deg \lambda_2 = 2(\ell +1)
\earr
\]
Hence we get a contradiction unless $(m_1, m_2)$ is $(1,0)$ or $(0,1)$.  The proposition now follows.
\hfill $\Box$

Proposition 3.1 will be used in Section 5 to determine the projective rank of the Grassmann manifolds.

\section{The Chen-Nagano Classification of H.S.S.}

In [5] B. Y. Chen and T. Nagano give a method of classifying compact symmetric spaces in terms of certain totally geodesic submanifolds.  In this section we summarize their method and in the next section apply it to the problem of determining the projective ranks of the compact irreducible Hermitian symmetric spaces.

Let $M$ be a symmetric space, $o  \in M$ the origin and $s_0$ the geodesic symmetry of $M$ fixing the point $o$.  A smooth closed geodesic through $o$, or circle $c$ for short, has an antipodal point $p$ with $s_o(p) =p$.  Write $M=G/K$ where $G$ is the group of isometries of $M$ (the closure of the group generated by all the symmetries $s_q, q \in M)$ and $K$ is the isotropy group of $o$.  Let $M_+(p)$ denote the orbit $K \cdot p$.  Then $M_+(p)$ is a totally geodesic submanifold of $M$ [5 2.1].  Let $\sigma_p$ be the involution of $G$ given by $\sigma_p(g) = s_pgs^{-1}_p$.  If $\gg = k + n$ is the Cartan decomposition of the Lie algebra $\gg =$ Lie $G$ with respect to $\sigma_0$ then $\sigma_p$ leaves $k$ stable ($s_0$ and $s_p$ commute) and induces the decomposition
\[	
	n=n_+ + n_-
\]
of $n$ into the $+1$ and $-1$ eigenspaces of $\sigma_p$.  Now $[[n_-, n_-],\  n_-] \subseteq n_-$ so $n_-$ forms a Lie triple system.  Consider the connected subgroup $G_-$ of $G$ corresponding to the Lie subalgebra $k_+ + n_-$.  Let $M_-(p)=G_- \cdot p$.  Then $M_-(p)$ is totally geodesic in $M$ and the tangent space to $M_-(p)$ at $p$ is the normal space to $M_+(p)$ at $p$ in $M$ (see [5, 2.2]).

Given a pair of antipodal points $(o,p)$ on a circle $c$ in $M$ we have the system $(o,p,M_+(p),M_-(p))$.  The group $G$ acts naturally on the set of all such systems.  Denote by $P(M)$ the orbit set.  If $f:B \rightarrow M$ is an isometric totally geodesic embedding then there arises a mapping $P(f): P(B) \rightarrow P(M)$ induced by the mapping carrying $(o,p,B_+(p), B_-(p))$ into $(f(o), f(p), M_+ (f(p)), M_-(f(p))$.  Moreover, $f(B_+(p)) \subset M_+(f(p))$ and $f(B_-(p)) \subseteq M_-(f(p))$.  We say $P(f)$ is a pairwise totally geodesic immersion.  The following theorem summarizes the results of Chen-Nagano which we shall use.

\begin{Thm} \quad [5, Section 5] Let $o$ be a point fixed by $K$ in the compact symmetric space $M=G/K$.  For $s \in G$ let $M^s$ denote the fixed point set of $s$.  Then
\end{Thm}

(1)  $M^{s_0} - \{0\}$ is the set of antipodal points on circles through $o$.

(2)	For each $p \in M^{s_0} = \{0 \}$ there exists an involutive automorphism $ad(b), \ b \in K \cap \exp(n)$ such that

\bigskip

\[
\barr{l}
	2(i) \quad M_+(p) = K(p) \ \mbox{is a covering space of} \ K/K^b \\
\noalign{\medskip}
	2(ii) \quad M_-(p), \ \mbox{the connected component of} \ M^b \ \mbox{containing} \ p  \\
	\mbox{is locally isometric with} \ G^b/K^b \\
\noalign{\medskip}
	2(iii) \quad \ \mbox{the tangent spaces to the totally geodesic submanifolds} \ M_+(p) \\
\noalign{\medskip}
	\mbox{and} \ M_-(p) \  \mbox{are the orthogonal complements of each other in} 
	\ T_p(M).
\earr
\]

\bigskip

(3)	The rank of $M_-(p)$ equals the rank $M$ and if $K$ is connected, the rank $K^b$ equals the rank of $K$.  

(4)	If a totally geodesic submanifold $B$ of $M$ has the same rank as $M$ then

\bigskip

\[
\barr{ll}
	(4(i)) & \qquad\qquad P(f) : P(B) \rightarrow P(M) \ \mbox{is surjection,} \ f \  \mbox{the inclusion} \\
\noalign{\medskip}
	(4(ii)) & \qquad\qquad \mbox{The Weyl group} \ W(B)\  \mbox{of}\ B \ \mbox{is a subgroup of} \ W(M)\\
\noalign{\medskip}
	(4(iii)) & \qquad\qquad \mbox{If} \ W(B) \  \mbox{is isomorphic to} \ W(M), \ \mbox{then}\ P(f) \ \mbox{is bijective}
\earr
\]

\bigskip

(5) $M$ is globally determined by $P(M)$; i.e., The set of isomorphism classes of compact irreducible symmetric spaces is in one-to-one correspondence with the set of the corresponding $P(M)$.

\bigskip

The following table gives the pairs $M_+(p)$ and $M_-(p)$ for the compact irreducible Hermitian symmetric spaces.  In this table $T,G^R(k,n-k), \ G^C(k,n-k)$, and $G^H(k,n-k)$ denote respectively the circle, the unoriented real Grassmann manifold $SU(n)/S(0(k) \times 0(n-k))$, the complex Grassmann manifold $SU(n)/S(U(k) \times U(n-k))$, and the quaternion Grassmann manifold $Sp(n)/Sp(k) \times Sp(n-k)$.  A determination is outlined in [5, p. 409].
\bigskip

\footnotesize

\noindent \begin{tabular}{llllc} \hline
	$M$ & $M_+$ & $M_-$ & Remark& $\#P(M)$ \\ \hline
\noalign{\medskip}
	$AIII(n):G^C(p,q)$ & $G^{\mathbb{C}}(h,p-h) \times G^{\mathbb{C}}(h,q-h)$ &
	$G^{\mathbb{C}}(h,h) \times  G^{\mathbb{C}}(p-h,q-h)$ & $0<h \leq p \leq q$ & $p$ \\
	$p+q=n$ \\
\noalign{\medskip}
	$BDI:G^{R}(p,q)$ & $G^R(h,p-h) \times G^{\mathbb{R}}(h,q - h)$ & 
	$G^R(h,h) \times G^R(p-h,q-h)$ & $0>h\leq p \leq q$ & $p$ \\
\noalign{\medskip}
	$CI(n):Sp(n)/U(n)$ & $G(h,n-h)$ & $CI(k) \times CI (n-k)$ & $0>k\leq n -1$ & $n$ \\
\noalign{\medskip}
	$DIII(n):SO(2n)/U(n)$ & $G^{\mathbb{C}}(k,n-k)$ & $DIH(k) \times DIH(n-k)$ & $0<k<n$ 
	& $[\frac{n}{2}]$ \\
\noalign{\medskip}
	$EIII$ & $DIII(5)$ & $S^2 \times G^C(5,1)$ & &$2$ \\
\noalign{\medskip}
	$EVII$ & $EIII$ & $S^2 \times G^R(10,2)$ & &$2$ \\  \hline 
\end{tabular}

\bigskip
\centerline{\bf Table 4.2}

\normalsize

\baselineskip20pt

\bigskip

\noindent{\bf 4.3 An Example} \quad Let $M = \bp^n_{\mathbb{C}}$.  Then $M = SU(n+1)/S[U(n) \times U(1)]$ and rank $(M) = 1$.  All geodesics in $M$ are circles, have the same length and are permuted transitively by $SU(n+1)$ [7, VII, 10.2, 11.2].  To find the possible geodesic pairs in $M$ it suffices to consider just one circle.  It is well-known that any projective line $L \subset \bp^n_{\mathbb{C}}$ (embedded linearly; i.e., of degree 1) is totally geodesic in $M$ so we may take the circle $c$ to lie in $L$.

Viewing $M$ as the set of $n$-dimensional subspaces of $\bc^{n+1}$, $L$ is the Schubert variety
\[
	L= \{W \in Gr(n,n+1) : W\supset W_0, \ W_0 \ \mbox{a fixed} \ n-1 \ \mbox{plane} \}
\]
Let $e_0, \ldots, e_n$ be an orthonormal basis for $\bc^{n+1}$ such that $W_0 = \langle e_0, \ldots, e_{n-2} \rangle$.  If $W \in L$ and $v$ is an element of $W$ then $v=w_0 + ae_{n-1} + be_n$ with $w_0 \in W_0$.  The assignment $v \longmapsto [a,b]$ where $[a,b]$ are homogeneous point coordinates in $\bp^1_{\mathbb{C}}$ defines an isomorphism of $L$ with $\bp^1_{\mathbb{C}}$.  We take our circle $c$ to be defined by
\[
	t \longmapsto \langle e_0, \ldots, e_{n-2}, \cos(t/2) \cdot e_{n-1} + \sin(t/2) \cdot e_n \}
\]
Let $g \in R = S[U(n) \times U(1)] \ \cap \ St_{SU(n+1)} \langle e_0, \ldots, e_{n-2}, e_n \rangle$.  Then since $W_0 = \langle e_0, \ldots, e_{n-2}, e_n \rangle \ \cap \ \langle e_0, \ldots, e_{n-2}, e_{n-1} \rangle, gW_0 = W_0$ so $R \subseteq H = St_{SU(n+1)} (W_0) \cap S[U(n) \times U(1)]$.  Hence we have a natural surjective map
\[
	Kp=K/R \rightarrow K/H \cong \bp^{n-1}_{\mathbb{C}}
\]
where $K=S[U(n) \times U(1)]$.  Since $K \cdot p$ is a connected totally geodesic subspace and is a proper subspace of $M$ and since $\bp^{n-1}_{\mathbb{C}}$ is simply connected this map is an isomorphism and $K \cdot p = \bp^{n-1}_{\mathbb{C}}$.  It follows that $(M_+(p), m_-(p)) = (\bp^{n-1}_{\mathbb{C}}, S^2 \simeq \bp^1_{\mathbb{C}})$ is the unique totally geodesic pair for $\bp^n_{\mathbb{C}}$.

\section{Projective Ranks of Compact H.S.S.}

We determine the projective ranks of the irreducible compact Hermitian symmetric spaces in this section.  We use the classification of Cartan as given in Table 4.2.

\subsection{Type AIII.}

We will show by induction on $n$ that $\mathbb{P}^{d+1}_{\mathbb{C}}$ embedded linearly in $Gr(d,n)$ is a maximal totally geodesic complex projective submanifold of $Gr(d,n)$.  If $n=1$ there is nothing to prove since $Gr(0,1) \cong \mathbb{P}^1_{\mathbb{C}}$.  Assume the result is known for all $Gr(d',n')$ with $n' < n$ where $n' \leq 2d'$ as usual.  Let $\mathbb{P}^K_{\mathbb{C}}$ be a maximal totally geodesic complex projective submanifold in $Gr(d,n)$.  Since $Gr(d,d+1) \simeq \mathbb{P}^{d+1}$ is totally geodesic in $Gr(d,n)$ when embedded linearly, $K \geq d+1$.  Suppose $K>d+1$.  Using Corollary 1.3 we can find a nontrivial homomorphism
\[
	\rho :SL(K +1, \mathbb{C}) \rightarrow SL(n+1, \mathbb{C})
\]
Since $2d \geq n, \ 2K \geq n$.  Thus $\mathbb{C}^{n+1}$ must split as an $SL(K + 1, \mathbb{C})$-module by Proposition 3.1.  Say $\cc^{n+1} \simeq \cc^{K+1} \oplus \cc^{n-k}$.  Let $0 \in Gr(d,n)$ correspond to the subspace $W \simeq \cc^{d+1} \subseteq \cc^{K+1}$.  Since $K > d+1, \ W$ has codimension at least $2$ in $\cc^{K+1}$.  Hence if $P$ is the stability group of $W$ in $SL(n+1, \cc)$ and $Q=P \cap SL(K+1, \cc)$ then $SL(K+1, \cc)/ Q \simeq Gr(d, K)$ is not a complex projective space because $K-d>2$.  This contradiction implies $K=d+1$.  Thus the projective rank of $Gr(d,n)$ is $d+1$.

Note that the maximal complex projective spaces $\pp^{d+1}_{\cc}$ in $Gr(d,n)$ are linearly embedded, are permuted transitively by $SU(n+1, \cc)$ and are parametrized by the $d+2$-dimensional subspaces of $\cc^{n+1}$, i.e., by $Gr(d+1,n)$.

\subsection{Type BDI.}

The Hermitian symmetric space
\[
	Q_m= SO(m+2)/SO(m) \times SO(2)
\]
is biholomorphic to the complex quadric in $\pp^{m+1}_\cc$
\[
	Q_m : z^2_0 + z^2_1 + \cdots + z^2_{m+1} = 0
\]
We have the following result from Chen-Nagano [4]:

\begin{SubThm} \quad Let $M$ be a totally geodesic complete connected Riemannian submanifold of $Q_m$.  Then 

\noindent 1. \ The embedding of $M$ in $Q_m$ is unique up to an isometry of $Q_m$. 

\noindent 2. \ If $M$ is maximal in $Q_m$ then $M$ is one of the following:

(i) \ $Q_m -1$

(ii) \ a local Riemannian product of two spheres $S^p \times S^q, \ p+q=m$

(iii) \ If $m=2n, \ \pp^n_{\cc}$ as Riemannian manifold

\noindent 3. \ If $M$ is not maximal then either $M \subseteq Q_{m-1}$ or $M \subseteq \pp^n(\rr)$ -- the real projective space of dimension $n$ if $m = 2n$. \hfill $\Box$
\end{SubThm}
The canonical decomposition of the Lie algebra $so(m + 2)$ is cf. [9, Vol. II, p. 278] $so(m+2) = \eta + p$ where
\[
	\eta = so(2) + so(m) = \left\{ \left(
	\barr{ccc}0 & -\lambda & 0 \\ \lambda & 0 & 0 \\ 0 & 0 & B \earr \right) : B \in so(m) \right.
\]
and
\[
	p= \left\{ \left( \barr{ccc} 0 & 0 & -t_\xi \\ 0 & 0 & -t_\eta \\ \xi & \eta & 0 \earr \right)
	: \xi, \eta \ \mbox{column vectors in} \  \rr^m \right\}
\] 
The complex structure of $p$ is given by $J(\xi, \eta)=(-\eta, \xi)$.  An inner product on $p \times p$ is given by
\[
	g(( \xi, \eta), (\xi', \eta')) = 4/c\{ \langle \xi, \xi' \rangle + 
	\langle \eta, \eta' \rangle \}
\]
where $\langle , \rangle$ denotes the standard Euclidean product of $\rr^n$ and $c$ is a positive scalar.  Both $g$ and $J$ are invariant under the adjoint action of $K = SO(2) \times SO(m)$ represented as
\[
	K= \left\{ \left[ \barr{ccc} \cos \theta & -\sin \theta & 0\cdots0 \\
	\sin \theta & \cos \theta & 0 \cdots 0 \\
	0 & 0 & \\
	\vdots & \vdots & B \\
	0 & 0 & \earr \right] : B \in SO(m) \right\}
\]
If $R(\theta) \times B$ is a typical element of $k$ as above and $(\xi, \eta) \in p = \rr^n \times \rr^m$ then
\[
	ad(\rr(\theta) \times B)(\xi, \eta) = (B \xi, B \eta)R(- \theta)
\]
where the right hand side is a matrix product.  The complex structure is the natural one so that
\[
	( \xi, \eta ) = \left( \barr{cc} \xi_1 & \eta_1 \\ \vdots & \vdots \\ 
	\xi_n & \eta_m \earr \right) \longmapsto \left( \barr{c} 
	\xi_1 + \sqrt{-1} \eta_1 \\ \vdots \\ \xi_n + \sqrt{-1} \eta _m \earr \right)
	\in \cc^n
\]
Thus the action of $SO(m) \subset K$ is the natural action of $SO(m)$ viewed as a real subgroup of $U(m)$.  The expression for $ad(R(\theta) \times B)( \xi, \eta)$ may be read as
\[
	ad(R(\theta) \times B)( \xi, \eta) = (B \cdot \underline{Z})e^{-i \theta}
\]
where $\underline{Z}$ is the complex vector with entries $\xi_i + \sqrt{-1}\eta_i$.

Let $m \subset p$ be a Lie triple system which is $J$-invariant.  Then $p=m \oplus m^\bot$ where the orthogonal splitting is relative to the standard Hermitian metric on $p$.  Since $g$ is the real part of this metric (up to scalar factor $4/c$) we have an exclusion of $SO(r), r = \dim m$, into $SO(m)$ and a Lie algebra homomorphism $so(r) \rightarrow so(m)$ compatible with adjoint actions on $m \subseteq p$.  This includes an inclusion of symmetric pairs
\[
	(so(r +2), so(r) \times so(2)) \rightarrow (so(m +2), so(m) \times so(2))
\]
and it follows readily that the totally geodesic submanifold determined by $m$ is the complex quadric $Q_r$ imbedded in $Q_m$ in some position.

\begin{SubProp} \quad The only totally geodesic complex submanifolds of $Q_m$ are the complex quadrics $Q_r \ r<m$ and the complex projective spaces $\pp^r_\cc, \ r \leq [\frac{m}{2}]$.
\end{SubProp}

\noindent{\bf Proof.} \quad Note that the assertion 2(iii) of 4.2 says that for $m = 2n, \ \pp^n_\cc$ is a totally geodesic submanifold.  When $m = 2n$, let $\rr^m = \rr^{2n}$ have a complex structure $J_0$.  Relative to this complex structure $SU (n)$ embeds in $SO(m)$.  Similarly for $\rr^{m+2} = \rr^{2n + 2} \simeq \cc^{n+1}$.  Any $n$-dimensional complex subspace has a natural induced orientation as a $2n$-dimensional real space.  Thus $G(n,n +1)$ can be mapped to $SO(m+2)/SO(2) \times SO(m)$ as follows: Let $e_1, \ldots, e_{n+1}$ be the standard orthonormal base for $\cc^{n+1}$ and for $W \subset \cc^{n+1}$ an $n$-plane write $\cc^{n+1} = W \oplus W^\bot$ as orthogonal direct sum.  Orient $W$ so that $W$ and $\langle e_1, \ldots, e_n \rangle$ have the same orientation -- $\langle e_1, \ldots, e_n \rangle$ and $W$ are in the same orbit under $SU(n+1)$.  Then orient $W^\bot$ so that $W \oplus W^\bot$ has the same orientation as $e_1, \ldots, e_{n+1}$.   Define the image of $W$ to be the oriented $2$-plane $W^\bot$ in $SO(m+2)/SO(2) \times SO(m)$.  The Lie algebra of $su(n+1)$ embeds in $so(m+2)$ by
\[
	A + iB \rightarrow \left( \barr{cc} A & B \\ -B & A \earr \right)
\]
from this one can check that $G(n,n+1) \rightarrow SO(m+2)/SO(2) \times SO(m)$ is a complex embedding.

\subsection{Type CI.}

Let $M = Sp(n)/U(n)$ be a H.S.S. of type CI.  Then $M$ is the subset of $G(n, 2n)$ consisting of $n$-dimensional subspaces of $\cc^{2n}$ which are totally isotropic with respect to a skew symmetric bilinear form [17, p. 232].  We take the matrix of this form to be 
\[
	J_n = \left[ \barr{cc} 0 & -J_n \\ J_n & 0 \earr \right]
\]
where $J_n$ is the $n \times n$-identity matrix.

Let $V_1 \subset \cc^{2n}$ be the subspace defined by $z_1 = \cdots = z_n = 0$ where $z_1, \cdots, z_{2n}$ are the usual coordinates of $\cc^n$.  Let $V_2 \subset \cc^{2n}$ be the subspace defined by the equations $z_i + z_{n+i} = 0 \quad 1 \leq i \leq n$.  Then $V_1$ and $V_2$ represent points in $M$ and $V_1 \cap V_2 = (0)$.  If we take the base point to be $p_0 = [V_1] \in M$ then $M = G_c / K$ where $G_c= Sp(n)/ ( \pm I_{2n})$ and $K=G_c \cap \{g \in Sp(n, \cc)| gp_0 = p_0\}$.  Then
\[
	K = \left\{ \pm \left( \barr{cc} A & 0 \\ 0 & \overline{A} \earr \right) : A \in U(n) \right\}
\]
It follows from this description that $M$ is a totally geodesic submanifold of $G(n, 2n)$.

According to Table 4.2 we expect that the projective rank of $M$ is $n-1$.  To see this consider $L \in G(n-1, V_1)$.  An easy calculation shows that $\dim(L^\bot \cap V_2) = 1$ and that $W = L \oplus (L^\bot \cap V_2)$ is totally isotropic in $\cc^{2n}$.  This induces a morphism $G(n-1, n, V_1) \rightarrow M$ exhibiting a map $\rho : \pp^{n-1}_\cc \rightarrow M$.  In fact, since $\cc^{2n} = V_1 \oplus V_2$ the above map is just the assignment $L \mapsto L + (L^\bot \cap V_2)$ for $L \subset V_1$ an $n-1$ plane.  Let $g_{12} \in Sp(n)$ be the element such that $g_{12}(V_1) = V_2$.  Then $SU(n) = SU(V_1)$ acts in $V_1 \oplus V_2$ by  $u(v_1, v_2) = (uv_1, g_{12} ug^{-1}_{12} v_2)$.  Under this action $SU(n)$ acts transitively on $\rho(G(n-1),V_1))$.  Since $\rho$ is $1-1$ and $SU(n)$ equivariant, $\rho(G(n-1, V))$ is totally geodesic in $M$.

Thus the projective rank of $CI(n)$ is at least $n-1$.  Since $M$ is totally geodesic in $G(n,2n)$ the maximal possible projective rank for $M$ is $n = pr(G(n,2n))$.  But if $B \subset \cc^{2n}$ is an $(n + 1)$-dimensional subspace then $G(n,B)$ cannot be contained in $M$.  Indeed, let $L \subset B$ be completely singular of codimension one and suppose $B=L + \cc_v$.  We have $L=L^\bot$ (this is easily seen for $L=V_1$ and since $Sp(n)$ acts transitively the same holds for $L$) choose codimension $1$ subspaces $L_1, \ldots, L_n$ of $L$ such that $L = L_1 + \cdots + L_n$ and put $L'_i = L_i + \langle v \rangle \subseteq B$.  If each $L'_i$ is completely singular then $v \in L^\bot_1 \cap L^\bot_2 \cap \cdots \cap L^\bot_n = (L_1 + \cdots + L_n)^\bot = L^\bot$ a contradiction.  Thus no $G(n,n +1)$ lies in $M$ and $pr(M) = n -1$.

\subsection{Type DIII}

Let $M = SO(2n)/U(n)$ be the H.S.S. of type DIII($n$).  Then $M$ may be identified with the submanifold if $G(n,2n)$ consisting of completely singular $n$-dimensional subspaces with respect to the form $S_n= \left[ \barr{cc} 0 & I_n \\ I_n & 0 \earr \right]$.  Let $e_1, \ldots, e_{2n}$ be the standard basis for $\cc^n$ and $V_1 = \langle e_1, \ldots, e_n \rangle$, $V_2 = \langle e_{n+1}, \ldots, e_{2n} \rangle$.  Then $V_1$ and $V_2$ are completely singular.  For $L \subset V_1$ a codimension one subspace, $\dim(L^\bot \cap V_2) = 1$ and $L + (L^\bot \cap V_2)$ is completely singular.  Thus we have a map $G(n-1, V_1) \rightarrow M$ which one checks as in 5.3 is totally geodesic.  Just as in 5.3, no $G(n,n+1)$ can lie in $M$ which is itself totally geodesic in $G(n, 2n)$.  Hence the projective rank of $M$ is also $n-1$.

\subsection{Type EIII}

The exception space is $M = E_6/SO(10) \cdot SO(2)$.  According to Table 4.2, $M$ contains the unique symmetric pair $M_+ = DIII(5), \ M_- = S^2 \times G^C(5,1)$.  We have seen that for any $K$ the symmetric pair for $\pp^K_{\cc}$ is $(\pp^{K-1}_{\cc}, S^2)$.  As the projective rank of DIII(5) is 4, by [5, p. 409] the pair $(\pp^4_\cc, S^2)$ is totally geodesically embeddable in $M$.  It follows that the maximal possible value of $pr(M)$ is 5.  From the Dynkin diagrams

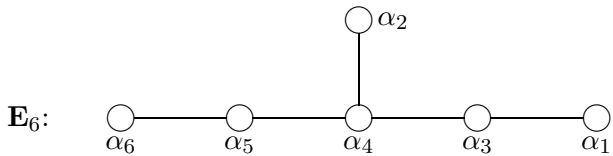
\begin{picture}(300,50)(-100,20)
\put (10,5){\circle{10}}
\put (10,-5){\makebox(0,0){$\alpha_6$}}
\put (15,5){\line(1,0){35}} 
\put (55,5){\circle{10}}
\put (55,-5){\makebox(0,0){$\alpha_5$}}
\put (60,5){\line(1,0){35}} 
\put (100,5){\circle{10}}
\put (100,-5){\makebox(0,0){$\alpha_4$}}
\put (100,42){\circle{10}}
\put (113,42){\makebox(0,0){$\alpha_2$}}
\put (100,10){\line(0,1){27}}
\put (105,5){\line(1,0){35}} 
\put (145,5){\circle{10}}
\put (145,-5){\makebox(0,0){$\alpha_3$}}
\put (150,5){\line(1,0){35}} 
\put (190,5){\circle{10}}
\put (190,-5){\makebox(0,0){$\alpha_1$}}
\put (-25,5){\makebox(0,0){{\bf E}$_6$:}}
\end{picture}

\null\vspace{-.60truein}
\begin{picture}(300,-25)(-100,40)
\put (10,5){\circle{10}}
\put (15,5){\line(1,0){35}} 
\put (55,5){\circle{10}}
\put (60,5){\line(1,0){35}} 
\put (100,5){\circle{10}}
\put (105,5){\line(1,0){35}} 
\put (145,5){\circle{10}}
\put (150,5){\line(1,0){35}} 
\put (190,5){\circle{10}}
\put (-25,5){\makebox(0,0){{\bf A}$_5$:}}
\end{picture}

\vspace{.75in}

We see that there exists a nontrivial homomorphism $s \ell(6, \cc)$ into $E_6$.  From [7, p. 507] the Cartan involution for $M$ is induced by the permutation (1,6) (3,5) of the vertices of the Dynkin diagram for $E_6$.  Thus $A_5$ is stable under the Cartan involution.  In particular, there exists a Lie triple system $m \subset E_6$ such that $m+[m,m]=s \ell(6, \cc)$.  This together with the fact from Table 4.2 that $G(5,1) \simeq \pp^5_\cc$ sits in $M$ tells us that indeed the projective rank of EIII is 5.

\subsection{Type E VII}

Let $M = E_7/E_6 \times SO(2)$ be the H.S.S. of type E VII.  The unique maximal symmetric pair for $M$ given by Table 4.2 is
\[
	(M_+, M_-) = (EIII, S^2 \times G^\rr(10,2))
\]
Since $pr(EIII) =5$ we have a symmetric pair
\[
	(\pp^5_\cc, S^2)
\]
contained in $E \ VII$ which suggest $pr(EV II) \leq 6$.

Again by considering the Dynkin diagram for $E_7$

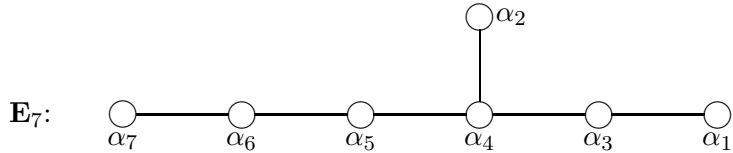
\begin{picture}(300,50)(-100,10)
\put (10,5){\circle{10}}
\put (10,-5){\makebox(0,0){$\alpha_7$}}
\put (15,5){\line(1,0){35}} 
\put (55,5){\circle{10}}
\put (55,-5){\makebox(0,0){$\alpha_6$}}
\put (60,5){\line(1,0){35}} 
\put (100,5){\circle{10}}
\put (100,-5){\makebox(0,0){$\alpha_5$}}
\put (145,42){\circle{10}}
\put (157,42){\makebox(0,0){$\alpha_2$}}
\put (145,10){\line(0,1){27}}
\put (105,5){\line(1,0){35}} 
\put (145,5){\circle{10}}
\put (145,-5){\makebox(0,0){$\alpha_4$}}
\put (150,5){\line(1,0){35}} 
\put (190,5){\circle{10}}
\put (190,-5){\makebox(0,0){$\alpha_3$}}
\put (195,5){\line(1,0){35}}
\put (235,5){\circle{10}}
\put (235,-5){\makebox(0,0){$\alpha_1$}}
\put (-25,5){\makebox(0,0){{\bf E}$_7$:}}
\end{picture}

\vspace{.5in}
\noindent we see that $A_6$ is in fact a subalgebra so again we have $s \ell (7, \cc) \subset E_7$.  Since $E_7 = k \oplus p$ with $k \simeq E_6 + \rr$ we have $(s \ell(7, \cc) \cap k)^c = s \ell(6, \cc) + \cc$.  Thus the corresponding totally geodesic subspace is of type AIII and is $G(6,7) \simeq \pp^6$.  Hence the projective rank of $E_7$ is $6$.  Summarizing we have the following result.

\begin{SubThm} \quad The projective ranks of the irreducible compact Hermitian symmetric spaces are as follows:
\begin{description}
\item{(i)} \quad pr[AIII(n,d)] = d, n/2 $\leq d < n$
\item{(ii)}\quad pr[BDI(m)] = [$\frac{m}{2}$]
\item{(iii)}\quad pr[CI(n)] = $n-1$
\item{(iv)} \quad pr[DIII(n)] = $n-1$
\item{(v)} \quad pr[EIII] = 5
\item{(vi)} \quad pr[EVII] = 6.
\end{description}
\end{SubThm}

\section{Degree and Conjugacy Results.}

Let $M=G/P$ be an irreducible Hermitian symmetric space with complex algebraic group of automorphisms $G$ and stability group of $0 \in M$ the group $P$.  As a $C$-space $M$ has a canonical ample line bundle $O_M(1)$ defined as follows:  If $P = L \cdot V$ is a Levi decomposition of $P$ with $V$ the unipotent radical of $P$, then $L$ is the connected centralizer in $P$ of a torus $S$.  Let $T$ be a maximal torus in $G$ containing $S$ such that the root system $\sum$ is relative to $T$.  If $\alpha$ is the root which defines the $C$-space, $S = (\mbox{ker} \  \alpha)^0$ and hence $\alpha$ induces  a character on $L$ -- the central character -- which will again be called $\alpha : L \rightarrow \cc^*$.  Let $E_\alpha = G \times_P \ \cc$ be the corresponding homogeneous line bundle on $M=G/P$.  Then $O_M(1)$ is the sheaf of sections of $E_\alpha$.  It is known that $O_M(1)$ is very ample and we measure the degree of a subvariety $Y \subset M$ relative to the  projective embedding defined by this line bundle.

Suppose $f : \pp^r_\cc \rightarrow M$ is a totally geodesic holomorphic embedding with image $Y$.  To find the degree of $f$, or $Y$, it suffices to find the degree of the restriction of $f$ to any projective line $L = \pp^1_\cc$ linearly embedded in $\pp^r_\cc$.  With a view toward the conjugacy problem in mind we begin with the following well-known result:

\begin{Prop} \quad There exists totally geodesic projective lines in $M$ of degree one.
\end{Prop}

\noindent{\bf Proof.} \quad We treat each type separately.  If $M$ is of type AIII the result follows from the fact that $G(d,d+1)$ has degree one in $Gr(d,n)$.

Suppose $M$ is of type BDI.  Viewing $M$ as the complex quadric $Q_n$ in $\pp^{n+1}_\cc$ we have a natural totally geodesic embedding of $Q_2$ in $Q_n$ as
\[
	Q_2 = \{ (Z_0, Z_1, Z_2, Z_3; 0, \ldots, 0) | Z^2_0 + Z^2_1 + Z^2_3 + Z^2_3 = 0 \}
\]
Consider the holomorphic isomorphism
\[
	\pp^1_\cc \times \pp^1_\cc \rightarrow Q_2
\]
given by
\[
	[u_0, u_1] \times [v_0, v_1] \stackrel{r}{\rightarrow} \left[\frac{u_0 v_0 + u_0 v_1}{2},
	\frac{u_0 v_1 - u_0 v_0}{2i}, \frac{u_1 v_0 - u_1 v_1}{2}, \frac{u_1 v_0 + u_1 v_1}{2i} \right]
\]
The image of the line $\pp^1 \times [0,1]$ is
\[
	L=\{ [u_0/2; u_0/2i; u_0/2i; u_1/2; u_1/2i ] \}
\]
which clearly has degree one in the ambient projective space.  By Mok [11, Section 1], $L$ is totally geodesic in $Q_2$ hence the proposition holds for $Q_n$ in general.

Now suppose $M$ is of type $CIII$ say $M = SP(n)/U(n)$.  As in Section 5 we view $M$ as $n$-dimensional completely singular subspaces of $\cc^{2n}$.  With $V_1$ and $V_2$ as in the last section we consider the mapping
\[
	\pp^1 \rightarrow G(n, 2n)
\]
given by $[u_0, u_1] \longmapsto \langle e_1, \ldots, e_{n-1}, u_0 e_n + u_1 e_{2n} \rangle$.  It is easily verified that the image lies in $M$.  Since the image also lies in $G(n,B)$ where $B= \langle e_1, \ldots, e_n, e_{2n} \rangle$ it readily follows that the degree of the image is one and, by Mok [{\textit ibid}] again, the image is totally geodesic in $M$.  When $M$ is of type DIII a completely analogous argument leads to the desired conclusion.

Finally, the exceptional cases EIII and EVII follow from the descriptions of the maximal totally geodesic complex projective spaces in these manifolds corresponding to the embeddings of $s \ell (6, \cc)$ and $s \ell(7, \cc)$ and the fact that the central characters in these cases restrict appropriately to $\mathcal{O}_{\pp^r_\cc}(1), r =$ projective rank of EIII (respectively EVIII). \hfill $\Box$

One might hope that the calculation of the degree of $f$ could be accomplished using the above proposition.  Unfortunately this is not the case.  Moreover, the degree of a totally geodesic embedding $\pp^s_\cc \rightarrow M$ is not a constant function of $s$.  The point is made with the following:

\medskip

\noindent{\bf Example 6.2} 
 \quad There exists a totally geodesic holomorphic embedding of $\pp^1_\cc$ into $Q_2$ of degree 2.

Consider the subalgebra $\sigma \simeq so(3)$ of $so(4)$ consisting of matrices
\[
	\sigma = \left\{ \left[ \barr{cc} A & 0 \\ 0 & 0 \earr \right] | A \in so(3) \right\}
\]
Then $\sigma$ is stable under the cartan involution $\theta$ of $Q_2$.  If $m(\sigma) \subset n$ is the $(-1)$-eigenspace then
\[
	m(\sigma) = \left\{ \left[ \barr{cccc} 0&0&-a&0 \\ 0&0&-b&0 \\ a&b&0&0 \\ 0&0&0&0 \earr
	\right] \right\}
\]
so $m(\sigma)$ is also $J$-stable.  The corresponding totally geodesic submanifold for this Lie triple system is therefore a complex submanifold.  Consider the two elements $X_1$ and $X_2$ of $\sigma$;
\[
	X_1 = \left[ \barr{cccc} 0&0&-1&0 \\ 0&0&0&0 \\ 1&0&0&0 \\ 0&0&0&0 \earr \right] \
	\mbox{and} \ X_2 = \left[ \barr{cccc} 0&0&0&0 \\ 0&0&-1&0 \\ 0&1&0&0 \\ 0&0&0&0 \earr \right]
\]
\hspace{1.5in} with $p = [1,i,0,0]$ as base point in $Q_2$

The geodesics given by $c_1 = \pi[\exp tX_1] \cdot p$ and $c_2 = \pi(\exp tX_2) \cdot p$ are
\[
	\barr{l} c_1 = \{ [1, i \cos t, i \sin t, 0 ] \; | \; t\in \rr \} \\
\noalign{\medskip}
	c_2= \{ [ \cos t, i, \sin t, 0 ] \; | \; t \in \rr \} \earr
\]
These curves lie on the algebraic curve
\[
	z_3 = z^2_0 + z^2_1 + z^2_2 = 0
\]
in $\pp^3_\cc$ which is the image (up to linear automorphism in $\pp^3_\cc$) of the veronese map of degree 2 mapping $\pp^1_\cc$ into the plane $z_3 = 0$ in $\pp^3_\cc$.  In particular, we have a totally geodesic complex projective line embedded in $Q_2$ having degree 2 and not 1. \hfill $\Box$

We turn now to the spaces of type CI.  Let $V_1$ and $V_2$ be the subspaces of $\cc^{2n}$ given in Section 5.3.  We have a morphism

\medskip

\noindent (6.3) \hspace{1.25in} $\rho : G(n -1, V_1) \times G(1, V_2) \rightarrow G(n,2n) = G(n, V_1+V_2)$

\medskip

\noindent given by $\rho(L, \ell) = L+ \ell$.  From the definition of the universal subbundle on $G(n, V_1 + V_2)$ is follows easily that $\rho^*(E(n, V_1 + V_2)) = p^*_1E(n-1, V_1) \oplus p^*_2 E(1, V_2)$ where $p_i$ denotes projection onto the $i$-th factor. 

Consider the morphism $f:G(n-1, V_1) \rightarrow G(1, V_2)$ given by $f(L) = L^\bot \cap V_2$.  Let $Z$ be the graph of this mapping in $G(n-1, V_1) \times G(1, V_2)$.  Then $\rho(Z)$ is precisely the totally geodesic submanifold of $M$ given in 5.3.  We want to compute the degree of the map $\tilde{f} = \rho \circ f_1$ where $f_1 : G(n-1, V_1) \simeq Z$.  More precisely we compute the determinant of the bundle $[\tilde{f}^* E(n, V_1 + V_2)]^\vee$.  Since $\rho^*E(n, V_1 + V_2) = p^*_1 E(n-1, V_1) \oplus p^*_2 E(1, V_2)$ it follows that $\tilde{f}^* E(n, V_1 + V_2) = f^*_1 p^*_1 E(n-1, V_1) \oplus f^*_1 p^*_2 E(1, V_2)$.  Evidently, $f^*_1 p^*_1 E(n-1, V_1) = E(n-1, V_1)$.  Now $E^\vee(1, V_2)$ is the line bundle corresponding to $\mathcal{O}(1)$ on $G(1, V_2) \simeq \pp^{n-1}_\cc$ and hence is generated by $n$-global sections.  Thus $f^*_1 p^*_2 E(1, V_2)^\vee$ is a line bundle on $G(n-1, V_1) \simeq \pp^{n-1}_\cc$ generated by $n$-global sections so $f^*_1 p^*_2 E(1, V_2) \simeq \mathcal{O} (-1)$.  It follows that $\tilde{f}^* E(n, V_1 + V_2) \simeq E(n-1, V_1) \oplus Q(n-1, V_1)^\vee$.  Hence
\[
	\barr{lcl}
	\det [\tilde{f}^* E(n, V_1 + V_2)]^\vee & = & \det[E(n-1, V_1)^\vee \oplus Q(n-1, V_1)] \\
	& = & \mathcal{O} (1) \oplus \mathcal{O}(1) \\
	& = & \mathcal{O}(2) 
	\earr
\]
An entirely analogous argument shows that if $M$ is of type DIII, then the degree of the embedding given in 5.4 is also 2.  Summarizing we have

\medskip

\noindent {\bf Proposition 6.4.} \quad \textit{Let $M$ be an irreducible compact Hermitian symmetric space of type CI or DIII.  Let $Y \subset M$ be a totally geodesic complex submanifold biholomorphic to $\pp^r_\pp, r = pr(M)$.  Then $\deg(Y) \geq 2$.}

\medskip

\noindent{\bf Proof.} \quad We have only to show that the case $\deg(Y) = 1$ cannot occur.  If $\deg(Y) = 1$ then by Theorem 2.3 there is an $(n+1)$-dimensional subspace $B \subset \cc^{2n}$ such that $Y$ is a hyperplane in $G(n,B)$.  Such a hyperplane consists of all $n$-dimensional subspaces of $B$ containing a fixed line, say $\cc \cdot v$.  But $v$ must be isotropic in this case.  As we have noted in 5.3, $B$ itself cannot be totally singular.  Thus there exists $w \in B$ with $J(w, w) \neq 0$.  The subspace spanned by $v$ and $w$ is contained in an $n$-dimensional subspace of $B$ so there exists at least one such $n$-plane which is not totally isotropic.  It follows that no hyperplane in $G(n,B)$ can be contained entirely in $M$ and hence $\deg(Y)$ cannot be one.  \hfill $\Box$

\medskip

Our next goal is to show that in fact the degree of $Y$ is precisely 2 when $M$ is of type $CI$ or DIII.  We will need a closer examination of totally geodesic complex projective subspaces of $G(n,2n)$ of dimension $n-1$.  Let $Y \subset G(n, 2n)$ be such a submanifold.  Then by 3.1 there exists a homomorphism $\rho: SU(n, \cc) \rightarrow SU(2n, \cc)$ and a corresponding map $SL(n, \cc) \rightarrow SL(2n, \cc)$.  Let $W_0 \subset \cc^{2n}$ be a fixed $n$-dimension subspace with orbit under $SU(n, \cc)$ isomorphic to $\pp^{n-1}_\cc$.  Let $H$ denote the isotropy group of $W_0$ in $SU(n, \cc)$.  Then $SU(n, \cc)/H \simeq \pp^{n-1}_\cc$.  

\noindent {\bf Lemma 6.5.} \quad \textit{Let $R$ be a compact subgroup of $SU(n, \cc)$ and suppose $SU(n, \cc)/R$ is biholomorphic to $\pp^{n-1}_\cc$.  Then $R$ is conjugate to $S[U(n-1) \times U(1)]$ in $SU(n, \cc)$.}

\medskip

\noindent{\bf Proof.} \quad Identify $\pp^{n-1}_\cc$ with $G(n-1, n)$ so that $S[U(n-1) \times U(1)]$ is the stability group of a point.  Then $R \subset S[U(n-1) \times U(1)]$.  By [35, Theorem 6.1ii)] R is a maximal connected proper subgroup of $SU(n, \cc)$ so the desired equality follows.  \hfill $\Box$

\medskip

The $SU(m, \cc )$-module $\cc^{2n}$ decomposes as a direct sum $V_1 \oplus V_2$ with $\dim V_1 = n$ and $V_2 = V^\bot_1$.  Since $H= S[U(n-1) \times U(1)] \simeq SU(n-1)$ is reductive each $V_i$ decomposes as a sum of irreducible $H$-modules:
\[
	V_i = W_i \oplus L_i, \ \dim W_i = n-1, \ \dim L_i =1
\]
Thus $\cc^{2n}$ decomposes as $H$-module
\[
	\cc^{2n} = W_1 \oplus L_1 \oplus W_2 \oplus L_2
\]

Now the submodule $W_j$ is $H$-stable and not fixed by $SU(n, \cc)$ so $W_0 \simeq W_i \oplus L_j \ i \neq j$ as $H$-modules.  More precisely we have an element $\tau \in H om_{H-mod} (W_i \oplus L_j, \cc^{2n})$ with image $W_0$.  Suppose for definiteness that $i = 1, \ j=2$.  Then
\[
	\barr{l}
	Hom_{H-mod} (W_i \oplus L_2, \cc^{2n}) \\
	\noalign{\medskip}
	= {\ds \prod_{i=1, 2}} Hom_{H-mod}(W_1 \oplus L_2, W_i \oplus L_1) \\
	\noalign{\medskip}
	\simeq {\ds \prod_{i=1,2}} [Hom_H (W_1, W_i) \times Hom_H(L_1, L_i)]
	\earr
\]
By Schur's lemma $\tau$ must be of the form $\tau(w, \ell) = (\beta w, 0, 0, \alpha, \ell)$ with $\alpha, \beta \in \cc^*$.  Thus $W_0 = W_1 \oplus L_2$ as $H$-module.

From this last equality we see that an element $g$ in $SU(n, \cc)$ maps $W_0$ into $gW_1 \oplus gL_2$.  It now follows that the orbit map $W_0 \rightarrow gW_0$ of $SU(n, \cc)$ into $G(n, 2n)$ is the same as the map described in 6.3 and hence $\deg(Y) = 2$ as desired. \hfill $\Box$

\medskip

We put the above calculations together and summarize the findings in the following results.

\medskip

\noindent {\bf Theorem 6.6.} \quad \textit{Let $M$ be an irreducible Hermitian symmetric space with connected isometry group $G_0$.  Let $Y$ be a totally geodesic complex projective subspace of $M,Y \simeq \pp^r_\cc$, with $r = pr(M)$.  Let $d$ be the degree of $Y$ in $M$.  Then
\begin{description}
\item{(i)} \quad $d=1$ if $M$ is of type AIII, EIII, or EVII.
\item{(ii)} \quad $d=2$ if $M$ is of CI or DIII.
\item{(iii)} \quad $d=1$ or $2$ if $M$ is of type BDI.
\end{description}}

\medskip

Moreover, all such submanifolds $Y$ of $M$ of minimal degree are conjugate under $G_0$.

\medskip

\noindent{\bf Proof.} \quad The assertions in (i) and (ii) have already been established.  As to (iii), by the discussion preceding 5.3 we see that each totally geodesic $\pp^1_\cc$ is contained in a complex quadric $Q_2 \simeq \pp^1_\cc \times \pp^1_\cc$.  Now $Q_2 \subseteq \pp^3_\cc$ as the hypersurface $z^2_0 + z^2_1 + z^2_3 + z^2_3 = 0$.  Since $\pp'_\cc$ is a complex submanifold or $Q_2$, its embedding in $\pp^3_\cc$ must be by a complete linear system $H^0(\pp^1, \mathcal{O}_\pp, (d))$.  If $d > 2$, then the forms of degree $d$ cannot satisfy the single quadratic relation above and yield an embedding so $d$ is at most 2 and (iii) follows.

As for the conjugacy assertion, again in the case of type AIII this follows immediately.  If $M$ is of type BDI, then we need only consider $Q_2 \subseteq M$.  Then the result follows from the fact that $SO(4, \rr)$ permutes the two-dimensional $J$-invariant Lie triple systems in $T_0, Q_2$ and clearly preserves the degree.

To treat the types CI and DIII we make use of the following:

\medskip

\noindent {\bf Lemma 6.7.} \quad \textit{Let $M$ be a compact irreducible Hermitian symmetric space of type AIII, CI or DIII and $j:\pp^1_\cc \rightarrow M$ a totally geodesic holomorphic isometric embedding of degree one.  Then every geodesic circle in $j(\pp^1_\cc)$ has minimal length.}

\medskip

\noindent{\bf Proof.}  \quad Since CI and DIII are totally geodesic in AIII it suffices to prove the result for $M =G(d,n)$.  For $G(d,n)$ the image of $j$ is contained in $G(d,d + 1)$ and a geodesic circle is given by (see example \#3)
\[
	t \longmapsto \langle e_0, \ldots, e_{d-2}, (\cos t/2) e_{d-1} + (\sin t/2)e_d \rangle
\]
$0 \leq t < 2 \pi$.  One checks easily that this has minimal length in $G(d,n)$ and the lemma follows.  \hfill $\Box$

\medskip

By [7, VII, 11.2] the closed geodesics of minimal length in $M$ are all conjugate.  Now in the case of CI($n$) or DIII($n$) we can consider the geodesic (in the notation of 5.3)
\[
	t \longmapsto \langle e_1, \ldots, e_{n-1}, ( \cos t/2) e_n + (\sin t/2)e_{2n} \rangle
\]
The mid-point corresponds to the subspace $\langle e_1, \ldots, e_{n-1}, ee_{2n} \rangle$ and the orbit under $SU(n, \cc)$ is precisely
\[
	\{ W \oplus L | W \subset V_1, L \subset V^\bot_1 \}
\]
i.e., the image of the map in 6.3.  Since the closed geodesics in $M$ are conjugate so are the submanifolds which arise from 6.3, 6.4.

Finally, we discuss the exceptional types EIII and EVII.  In each case according to Table 4.1 there exists a unique Hermitian pair $(M_+, M_-)$.  Considering first the case EIII,
\[
	(M_+, M_-) = (DIII(5), S^2 \times G(1,6))
\]
If $Y_1, Y_2$ are the images of two totally geodesic embeddings of $\pp^5_\cc$ in $M$, then we may first assume $(Y_{1+}, Y_{1-} ) = (Y_{2+}, Y_{2-})$ by the result for DIII.  Then $Y_1$ and $Y_2$ meet along a hyperplane and have the same normal space in EIII so coincide.

If $M = EVII$, then
\[
	(M_+ , M_-) = (EIII, S^2 \times G^\rr(12,2))
\]
Again for $Y_1, Y_2$ images of totally geodesic maps from $\pp^6_\cc$ into $M$ we may assume $(Y_{1+}, Y_{1-}) = (Y_{2+}, Y_{2-})$ so $Y_1$ and $Y_2$ meet along a $\pp^5_\cc$ and have the same normal space in EVII so coincide.  This completes the proof.

\pagebreak

\end{document}